\documentclass[12pt]{article}
\usepackage{graphicx}
\usepackage{amssymb}

\begin{document}

\begin{center}
{\bf Counting Triangles in Triangles} \\
\ \\
Jim Propp (UMass Lowell, Lowell, MA) \\
Adam Propp-Gubin (Belmont High School, Belmont, MA)
\end{center}

\begin{abstract}
\noindent
We give a formula for counting the triangles in a picture 
consisting of the three sides of a triangle and some cevians. 
This lets us prove statements that are claimed without proof in 
the Online Encyclopedia of Integer Sequences and some popular YouTube videos,
and also prove some new results. We also give formulas that apply 
when the cevians cut each side into equal-length pieces.
\end{abstract}

\section{Introduction}

How many triangles can you find in Figure 1?
\begin{figure}[hp]
\begin{center}
\includegraphics[width=1.5in]{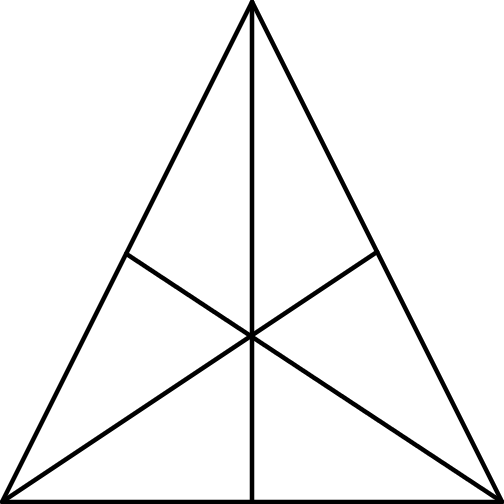}
\caption{A triangle with its three medians.}
\end{center}
\end{figure}
Keep in mind that triangles can be made of smaller triangles,
like the two shown in Figure 2.
\begin{figure}[h]
\begin{center}
\includegraphics[width=1.5in]{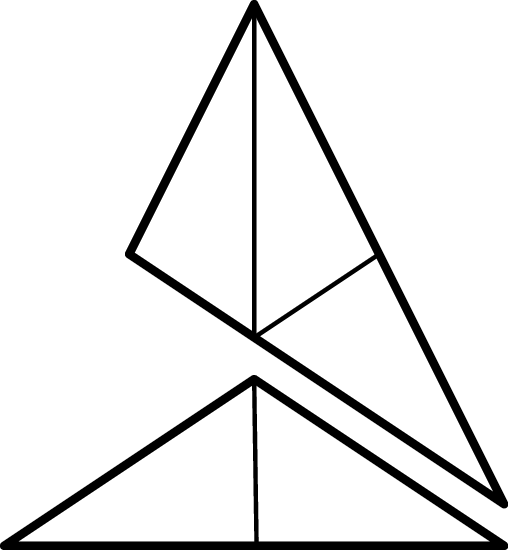}
\end{center}
\caption{Two triangles in a triangle.}
\end{figure}
Do you think you found all the triangles in the original picture? 
Are you sure?
Did you forget to include the big triangle that contains
all six of the little triangles?
Now that you've included it,
do you think you've found all the triangles in the picture?
How can you be sure?

There are a lot of YouTube videos that treat such problems, 
such as a Presh Talwalkar's 
``How Many Triangles Are There?'' video~\cite{T}.
Talwalkar's video, which has over half a million views,
is about counting triangles in Figure 3.
\begin{figure}
\begin{center}
\includegraphics[width=1.5in]{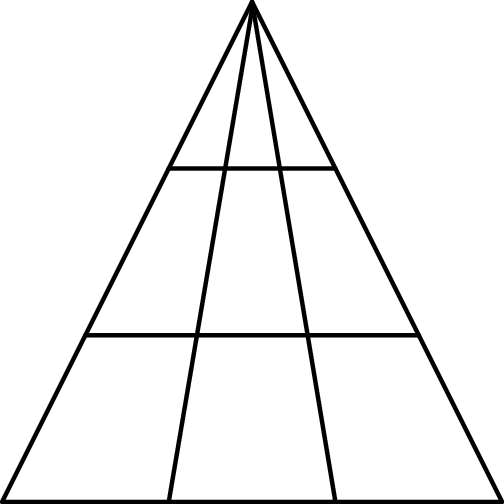}
\caption{The ``Bollywood puzzle''.}
\end{center}
\end{figure}
Problems like this appear on a lot of Indian examinations.
Talwalkar's video is better than most videos of its kind
because it uses mathematical reasoning.

A more typical video is Imran Sir's 
``Best Trick for Counting Figures'' video~\cite{S}, 
which gives tricks for getting the right answer quickly
without any explanation for why the tricks work.
For instance, Sir claims that 
the number of triangles in Figure 4 is $4^3$.
\begin{figure}
\begin{center}
\includegraphics[width=1.5in]{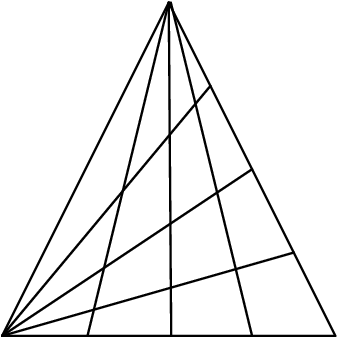}
\end{center}
\caption{Sir's puzzle.}
\end{figure}
The six extra lines inside the triangle are {\em cevians}:
lines that pass through a vertex of a triangle 
and a point on the opposite side.
Sir makes the more general claim that if you have
a triangle with $n$ cevians coming from one vertex
and $n$ cevians coming from a second vertex
(and no cevians coming from the third vertex)
then the number of triangles in the picture is $(n+1)^3$.
But why? This claim also appears
in the entry for the sequence of perfect cubes 
in the Online Encyclopedia of Integer Sequences~\cite{O1},
where it's credited to Lekraj Beedassy,
but as far as we know no proof has been published.
The OEIS also includes a formula for counting triangles
in a more complicated picture in which there are $n$ cevians coming 
from each of the three vertices of a triangle
on the assumption that no three cevians meet at a point~\cite{O2},
but the argument given there isn't very clear.

The purpose of this article is to prove all these claims
by providing a general formula for counting triangles
in a picture consisting of a triangle and some cevians.
The proofs are easy but we feel that
because of the global interest in such puzzles,
there should be a published article that explains what's going on
in this special class of triangle-counting puzzles.
Our approach also provides a nice application of 
the combinatorial method of overcounting-and-correcting.

Our main result is 

\bigskip

{\bf Theorem 1}: In triangle $ABC$, 
if you draw $a$ cevians from vertex $A$,
$b$ cevians from vertex $B$, and $c$ cevians from vertex $C$,
then the number of triangles formed is
$${a+b+c+3 \choose 3}-{a+2 \choose 3}-{b+2 \choose 3}-{c+2 \choose 3}-d$$
where $d$ is the number of points inside the triangle
where three cevians meet.

\bigskip

Here ${n \choose 3}$ is the number of ways to choose 3 objects
out of a collection of $n$ objects; it equals $\frac{n!}{3!(n-3)!}$,
which simplifies to $n(n-1)(n-2)/6$.

Consider again the problem that we gave at the beginning of this article.
We could draw the sixteen pictures shown in Figure 5,
and after checking that there are no duplicates,
we'd know that the answer is at least sixteen.
\begin{figure}[t]
\begin{center}
\includegraphics[width=4in]{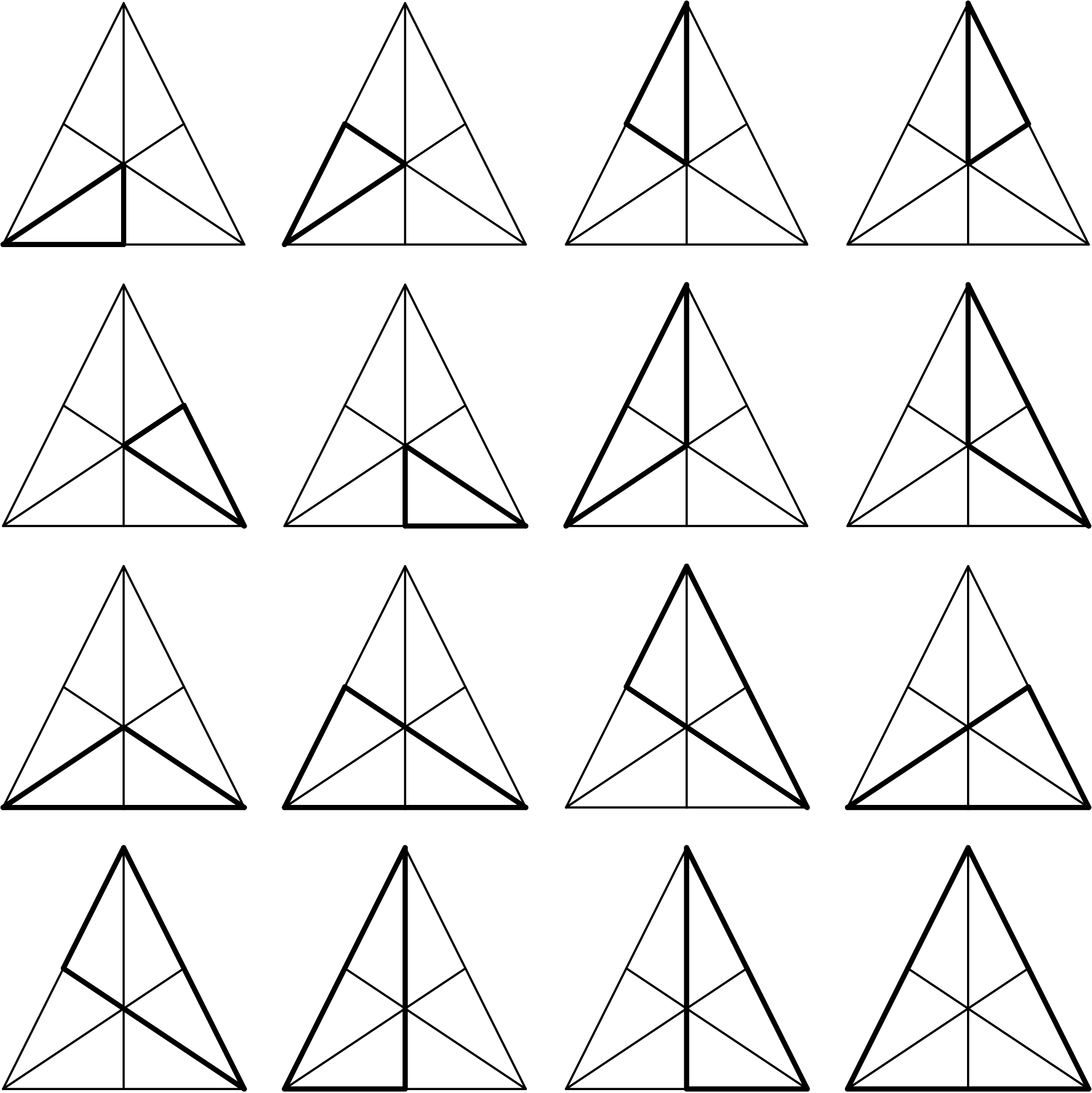}
\end{center}
\caption{Sixteen triangles in a triangle.}
\end{figure}
But how would we know we didn't miss any other triangles?
Theorem 1 assures us that we didn't:
we have $a=b=c=1$ and we can see that $d=1$ as well, so that the formula gives
${6 \choose 3} - {3 \choose 3} - {3 \choose 3} - {3 \choose 3} - 1 = 16$.

Notice that $d$ must be 0 when $a=0$ or $b=0$ or $c=0$,
since the only way three cevians can meet inside the triangle
is if one cevian is an $A$-cevian, one is a $B$-cevian, and one is a $C$-cevian.
Because of this, when $a=b=n$ and $c=0$, we automatically get $d=0$
so that the formula gives
${2n+3 \choose 3} - {n+2 \choose 3} - {n+2 \choose 3} - {2 \choose 3} - 0$,
which simplifies to $(n+1)^3$, as Sir and Beedassy claim.
But there's no need for $a$ and $b$ to be equal if $c=0$;
we still automatically get $d=0$,
so Theorem 1 tells us that the number of triangles is
${a+b+3 \choose 3}-{a+2 \choose 3}-{b+2 \choose 3}$,
which simplifies to $(a+1)(b+1)(a+b+2)/2$.

Another application of Theorem 1 is the case $a=b=n$, $c=1$
with the assumption that the picture has bilateral symmetry
under the reflection that fixes $C$ and switches $A$ and $B$.
Figure 6 shows the case $n=3$.
\begin{figure}
\begin{center}
\includegraphics[width=1.8in]{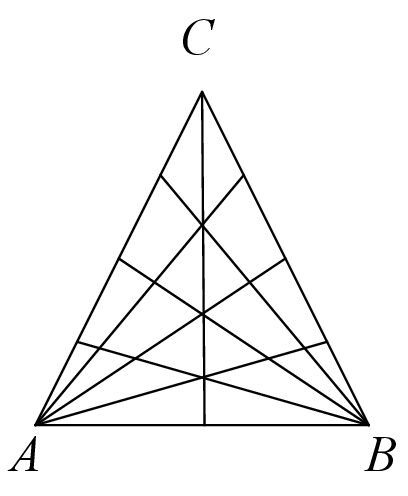}
\end{center}
\caption{A puzzle with bilateral symmetry.}
\end{figure}
In this situation there are $n$ points where three cevians meet,
corresponding to triples that consist of
a cevian from $A$ to $BC$, its mirror-twin (from $B$ to $AC$),
and the median from $C$ to $AB$. Our formula gives
${2n+4 \choose 3} - 2 {n+2 \choose 3} - {3 \choose 3} - n
= (n+3)(n+1)^2$; this provides a new combinatorial
interpretation of this sequence~\cite{O3}.

Theorem 1 bypasses the sometimes difficult problem of finding $d$,
which can't be determined from $a$, $b$, and $c$ alone.
For instance, in Figure 7 we have $a=b=c=1$
just like in the opening example, but $d$ is 0 instead of 1,
so that the number of triangles is 17 instead of 16.
\begin{figure}
\begin{center}
\includegraphics[width=1.5in]{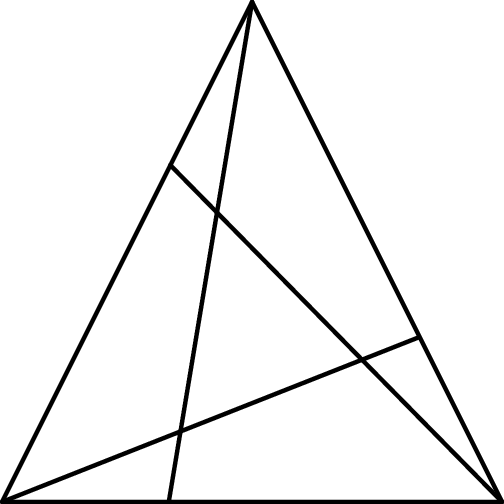}
\end{center}
\caption{Three cevians that don't intersect.}
\end{figure}

Theorem 2 establishes an important case in which 
the geometry of the cevians is specified
in enough detail to allow us to
determine $d$ and so find a formula
for the number of triangles in the picture:

\bigskip

{\bf Theorem 2}: Let $p$ be a prime
and let $q=p^m$ for some exponent $m \geq 1$.
From each vertex of triangle $ABC$ draw $q-1$ cevians
that divide the opposite side into $q$ pieces of equal length.
Then 

(a) If $p$ is odd, the number of triangles formed is 
$(8q^3 - 9q^2 + 3q)/2$.

(b) If $p=2$, the number of triangles formed is 
$(8q^3 - 9q^2 - 3q + 10)/2$.

\bigskip

Our proof of Theorem 2 relies on Theorem 1
along with Ceva's Theorem from geometry
(from which cevians got their name),
Euclid's Lemma from number theory (when $q$ is an odd prime),
and the Fundamental Theorem of Arithmetic 
(when $q$ is a general prime power).

\section{Proofs}

{\bf Proof of Theorem 1}:
There are $a+b+c+3$ line segments in total
(the cevians along with the three sides of triangle $ABC$).
No two of the segments are parallel,
and each of the segments intersects each of the others,
either in the interior of triangle $ABC$
or else at one of the three corners.
Each of the triangles in the picture 
is bounded by three of these $a+b+c+3$ line segments,
and each triple of line segments either bounds a triangle
or else meets at a single point.
So the number of triangles can't be more than ${a+b+c+3 \choose 3}$;
to turn this into the right answer, we have to subtract 
the number of triples of line segments that meet at a point. 
There are $a+2$ line segments that pass through $A$,
so there are $a+2 \choose 3$ triples of line segments that meet at $A$.
For the same reason 
there are $b+2 \choose 3$ triples of line segments that meet at $B$
and $c+2 \choose 3$ triples of line segments that meet at $C$.
Lastly, there are $d$ triples of line segments that meet in the interior,
since each such intersection point
uniquely determines the three cevians passing through it.
So the number of triples of line segments that bound a triangle is
${a+b+c+3 \choose 3} - {a+2 \choose 3} - {b+2 \choose 3} - {c+2 \choose 3} - d$,
as claimed.
$\square$

\bigskip

\noindent
{\bf Proof of Theorem 2}:
By Ceva's Theorem, in any triangle $ABC$ (see Figure 8),
\begin{figure}
\begin{center}
\includegraphics[width=3in]{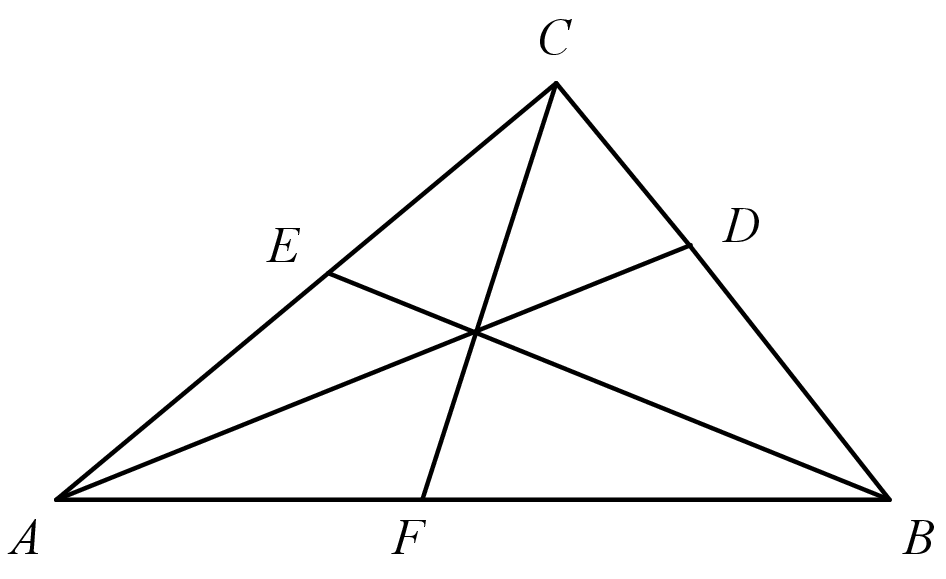}
\caption{Ceva's Theorem.}
\end{center}
\end{figure}
the three cevians $AD$, $BE$, and $CF$ meet at a point
if and only if $$\frac{BD}{DC} \frac{CE}{EA} \frac{AF}{FB} = 1.$$
The $p-1$ cevians from $A$ divide segment $BC$ into $p$ equal parts,
so if $AD$ is one of those cevians,
line segment $BD$ must consist of some whole number of those equal parts, 
say $i$ of them, so that the complementary line segment
$DC$ must consist of $p-i$ of them.
So the ratio $BD/DC$ must be of the form $i/(p-i)$
for some integer $i$ between 1 and $p-1$. 
For the same reason there must be an integer $j$ between 1 and $p-1$
for which $CE/EA = j/(p-j)$ 
and an integer $k$ between 1 and $p-1$
for which $AF/FB = k/(p-k)$.
So the three cevians associated with $i$, $j$, and $k$
meet at a point if and only if $(i/(p-i)) (j/(p-j)) (k/(p-k)) = 1$,
that is, if and only if $i j k = (p-i) (p-j) (p-k)$.
Expanding the right hand side and moving $ijk$ to the left side,
we can rewrite this as $$2ijk = p^3 - p^2 (i+j+k) + p (ij+ik+jk).$$

At this point the analysis of the equation splits into two cases.

(a) Assume $p$ is odd. First we do the sub-case where $m$ is 1, 
that is, $q$ is equal to an odd prime $p$.
Since $p$ divides the right hand side of the preceding equation,
$p$ must divide the left hand side as well.
But since $p$ is an odd prime, $p$ can't divide 2,
nor can $p$ divide $i$, $j$, or $k$
since $p$ is bigger than all three of them.
On the other hand, Euclid's lemma tells us that
when a prime $p$ divides a product of two or more numbers,
it must divide at least one of them.
This contradiction shows that no triple $i$, $j$, $k$ 
with $1 \leq i,j,k \leq p-1$
satisfies $i j k = (p-i) (p-j) (p-k)$.
That is, there are no points interior to the triangle
where three cevians intersect. So $d$ is 0 as claimed,
and Theorem 1 (with $a=b=c=p-1$)
then tells us that the number of triangles is
${3p \choose 3} - 3 {p+1 \choose 3}$,
which simplifies to $(8 p^3 - 9 p^2 + 3p)/2$,
as our formula predicts.

In the more general situation where $q=p^m$ (with $p$ odd)
we apply the Fundamental Theorem of Arithmetic,
which says that every natural number 
can be written in a unique way as a product of powers of primes
if we don't pay attention to the order
in which the powers of primes get multiplied together.
Because of this, every natural number $n$ can be written in a unique way
as a power of $p$ times a number not divisible by $p$.
Write $i = p^r f$, $j = p^s g$, and $k = p^t h$
where $f$, $g$, and $h$ are integers not divisible by $p$.
It's important to notice that $r$, $s$, and $t$ 
must all be strictly less than $m$ 
(this is where the assumption that $q$ is a power of $p$ plays a role).
We can rewrite the equation $i j k = (q - i) (q - j) (q - k)$
as $2 i j k = q^3 - q^2 i - q^2 j - q^2 k + q i j + q i k + q j k$
and then rewrite it again as
$$2 p^{r+s+t} e f g = p^{3m} - p^{2m+r} f - p^{2m+s} g - p^{2m+t} h
+ p^{m+r+s} fg + p^{m+r+t} fh + p^{m+s+t} g h.$$
At first it's not clear how this will help us,
but all the exponents of $p$ in the right hand side
are bigger than $r+s+t$ because $m$ is greater than $r$, $s$, and $t$.
(For instance, if we add together the three inequalities
$m>r$, $m>s$, and $m>t$ we get the inequality $3m>r+s+t$,
which tells us that the first exponent in the right hand side
is bigger than $r+s+t$;
or, if we add $s+t$ to both sides of the inequality $m > r$ 
we get the inequality $m+s+t > r+s+t$
which tells us that the last exponent in the right hand side
is bigger than $r+s+t$.)
So each exponent in the right hand side
is strictly bigger than $r+s+t$,
or in other words, each of those exponents
is greater than or equal to $r+s+t+1$.
So each term in the right hand side of the equation
is divisible by $p^{r+s+t+1}$, which implies that the whole
right hand side is divisible by $p^{r+s+t+1}$.
However, the left hand side isn't divisible by $p^{r+s+t+1}$
because the exponent in $p^{r+s+t}$ falls short by 1
and because 2, $e$, $f$, and $g$ are not divisible by $p$
(here we're using Euclid's lemma again).
This contradiction tells us that $d=0$ in this case as well,
so the formula we want to prove follows by the same reasoning
that we used in the proof of the case where $q$ is prime.

(b) Now assume that $p = 2$, so that $q = 2^m$.
The equation $ijk = (q-i)(q-j)(q-k)$ becomes 
$$2^{r+s+t+1} e f g = 2^{3m} - 2^{2m+r} f - 2^{2m+s} g - 2^{2m+t} h
+ 2^{m+r+s} fg + 2^{m+r+t} fh + 2^{m+s+t} g h.$$
Because $r$, $s$, and $t$ are all at most $m-1$,
the exponent $3m$ must be at least 3 more than $r+s+t$
and the exponents $2m+r$, $2m+s$, and $2m+t$
must all be at least 2 more than $r+s+t$,
so none of those exponents can equal $r+s+t+1$.
The only exponents that might equal $r+s+t+1$ are the last three:
$m+r+s$, $m+r+t$, and $m+s+t$.
Suppose without loss of generality that
$m+s+t$ (the last of them) equals $r+s+t+1$.
Then we have $m=r+1$, or in other words $r=m-1$.
So $i$ must be a multiple of $2^{m-1}$ that's smaller than $2^m$,
and the only possibility is $i = 2^{m-1}$ itself (that is, $i = q/2$).
In other words, the $i$th cevian from vertex $A$ is actually a median.
This shows that at every point where cevians intersect, 
at least one of the three cevians must be a median.
But there are three medians to consider.
For each median there are $q-1$ mirror-pairs of cevians,
so we might think the total number of intersection points
is three times that; but the expression $3(q-1)$
triple-counts the point where all three medians meet, so we have to subtract 2.
So the total number of intersection points is $3(q-1)-2 = 3 q - 5$.
Plugging this into Theorem 1, we find that the number of triangles is
$(8 q^3 - 9 q^2 - 3 q + 10)/2$.
$\square$

\section{Conclusion}

Theorem 1 was proved using the strategy
of overcounting followed by correcting.
This strategy can also be applied to the ``Bollywood problem'' 
mentioned near the start of the article.
We've got 7 lines, so there are ${7 \choose 3} = 35$ triples.
${4 \choose 3} = 4$ of the triples meet at the apex,
${3 \choose 3} = 1$ of the triples consist of three parallel line segments,
and ${3 \choose 2} {4 \choose 1} = 12$ of the triples
consist of two parallel line segments and a line that goes through the apex.
So the number of triangles in this picture is $35-4-1-12$, or 18.
A more direct approach like Presh Talwalkar's yields the answer with less work,
but it's always good to know multiple ways to arrive at
the answer to a problem, especially on a math test,
if you have enough time to check your work. 

One problem that proved to be too hard for us
was finding an exact formula that's like Theorem 2
in assuming that the cevians divide the opposite side
into equal-length pieces but doesn't assume that 
the number of cevians is 1 less than a power of a prime.
The difficult part is counting the places where cevians intersect.
That is, for general values of $n$
it's hard to determine the number of solutions to the equation
$ijk = (n-i)(n-j)(n-k)$ with $i,j,k$ between 1 and $n-1$.
This mysterious number, which depends on the prime factorization of $n$
in a complicated way, is the subject of sequence A331423 in the OEIS~\cite{O4}.
A related entry is sequence A332378\cite{O5};
it lists the odd values of $n$
for which there's a triple of cevians meeting at a point
(under the assumption that the cevians from each vertex of the triangle
divide the opposite side into $n$ pieces of equal length).

We did notice two patterns in the data.
If $q$ is of the form $p(2p-1)$ where $p$ and $2p-1$ are both prime,
then it seems that the picture you get from having $3(q-1)$ cevians 
dividing each side of the triangle into $q$ pieces of equal length
gives a positive number of intersection points;
we've checked this out to $p=97$, $p(2p-1)=18721$.
The same sort of thing seems to happen
when $q$ is of the form $p^2(2p+1)$ where $p$ and $2p+1$ are both prime
(that is, where $p$ is a Germaine prime); we've checked this out to 
$p=29$, $p^2(2p+1)=49619$.
% Still running on Monday morning!
We don't know how to prove that these patterns continue, let alone 
find a general rule governing those two integer sequences.
Maybe you can!

\section*{Acknowledgement}

We thank Richard Green for pointing out a mistake in this 
article: the picture of sixteen triangles-in-a-triangle
omits one of the possibilities and contains another
possibility twice. The first author takes full responsibility 
for this mistake. Rather than correct the error,
we have chosen to leave it as-is, as a vivid demonstration
of the pitfalls of counting-by-listing. (Did you notice the mistake?)

\bigskip

\noindent
{\sc Jim Propp} is a mathematician at UMass Lowell
with interests in combinatorics, probability, and
writing about mathematics for the public. He is the
father of Adam Propp-Gubin.

\medskip

\noindent
{\sc Adam Propp-Gubin} is a senior at Belmont High School.
He saw some videos about triangle-counting problems and
asked his father ``Is this the sort of thing you can do
math research on?'' This article is the result..

\medskip
\end{document}